\documentclass[11pt,a4paper]{article}
\usepackage{amsfonts,amsmath,amssymb,accents,amsthm}
\usepackage{verbatim}

\usepackage{hyperref}
\usepackage[normalem]{ulem}

\usepackage{tikz} 

\newcommand{\dis}{\displaystyle}

\newcommand{\noi}{\noindent}
\newcommand{\halmos}{\rule{1ex}{1.4ex}}
\newcommand{\QED}{\nopagebreak{\hspace*{\fill}$\halmos$\medskip}}

\newcommand{\med}{\medskip}
\newcommand{\quand}{\quad\mbox{and}\quad}

\newtheoremstyle{mythm}
  {}
  {}
  {\itshape}
  {}
  {\bfseries}
  {}
  {.5em}
  {#1 #2 \thmnote{(#3)}}

\theoremstyle{mythm}
\newtheorem{theorem}{Theorem}
\newtheorem{proposition}[theorem]{Proposition}
\newtheorem{lemma}[theorem]{Lemma}
\newtheorem{exercise}[theorem]{Exercise}
\newtheorem{corollary}[theorem]{Corollary}
\newtheorem{conjecture}[theorem]{Conjecture}

\newtheorem{counterex}[theorem]{Counterexample}

\newcommand{\bt}{\begin{theorem}}
\newcommand{\et}{\end{theorem}}
\newcommand{\bl}{\begin{lemma}}
\newcommand{\el}{\end{lemma}}
\newcommand{\bp}{\begin{proposition}}
\newcommand{\ep}{\end{proposition}}
\newcommand{\bcor}{\begin{corollary}}
\newcommand{\ecor}{\end{corollary}}
\newcommand{\br}{\begin{remark}\rm}
\newcommand{\er}{\end{remark}}
\newcommand{\bcon}{\begin{conjecture}}
\newcommand{\econ}{\end{conjecture}}
\newcommand{\bex}{\begin{exercise}}
\newcommand{\eex}{\end{exercise}}
\newcommand{\bcou}{\begin{counterex}}
\newcommand{\ecou}{\end{counterex}}

\newenvironment{Proof}[1][]{\noi\textbf{Proof #1}}{\QED}
\newcommand{\bpro}{\begin{Proof}}
\newcommand{\epro}{\end{Proof}}

\newcommand{\be}{\begin{equation}}
\newcommand{\ee}{\end{equation}}
\newcommand{\ba}{\begin{array}}
\newcommand{\ea}{\end{array}}
\newcommand{\bc}{\be\begin{array}{r@{\,}c@{\,}l}}
\newcommand{\ec}{\end{array}\ee}

\newcommand{\ga}{\gamma}

\newcommand{\de}{\delta}

\newcommand{\eps}{\varepsilon}
\newcommand{\la}{\lambda}

\newcommand{\sig}{\sigma}

\newcommand{\om}{\omega}
\newcommand{\Om}{\Omega}

\newcommand{\Ai}{{\cal A}}

\newcommand{\Di}{{\cal D}}

\newcommand{\Pc}{{\cal P}}

\newcommand{\Ui}{{\cal U}}

\newcommand{\R}{{\mathbb R}}

\newcommand{\E}{{\mathbb E}}
\renewcommand{\P}{{\mathbb P}}

\newcommand{\li}{\langle}
\newcommand{\re}{\rangle}

\newcommand{\desd}{\ensuremath{\Leftrightarrow}}

\newcommand{\up}{\uparrow}
\newcommand{\down}{\downarrow}
\newcommand{\sub}{\subset}
\newcommand{\beh}{\backslash}

\newcommand{\ffrac}[2]{{\textstyle\frac{{#1}}{{#2}}}}
\newcommand{\dif}[1]{\ffrac{\partial}{\partial{#1}}}
\newcommand{\diff}[1]{\ffrac{\partial^2}{{\partial{#1}}^2}}

\newcommand{\di}{\mathrm{d}}
\newcommand{\half}{{[0,\infty)}}
\newcommand{\expo}{\mbox{\large\it e}}
\newcommand{\ex}[1]{\expo^{\,\textstyle{#1}}}
\newcommand{\ha}{\ffrac{1}{2}}

\newcommand{\var}{{\rm Var}}

\newcommand{\hatOm}{\widehat\Om}

\begin{document}

\makeatletter\@addtoreset{equation}{section}
\makeatother\def\theequation{\thesection.\arabic{equation}} 

\renewcommand{\labelenumi}{{\rm (\roman{enumi})}}
\renewcommand{\theenumi}{\roman{enumi}}

\title{Necessary and sufficient conditions for a nonnegative\\ matrix to be
  strongly R-positive}
\author{Jan~M.~Swart\footnote{The Czech Academy of Sciences,
Institute of Information Theory and Automation,
Pod vod\'arenskou v\v e\v z\' i 4,
18200 Praha 8,
Czech Republic;
swart@utia.cas.cz}
}

\date{\today}

\maketitle

\begin{abstract}\noi
Using the Perron-Frobenius eigenfunction and eigenvalue, each finite
irreducible nonnegative matrix $A$ can be transformed into a probability
kernel $P$. This was generalized by David Vere-Jones who gave necessary and
sufficient conditions for a countably infinite irreducible nonnegative matrix
$A$ to be transformable into a recurrent probability kernel $P$, and showed
uniqueness of $P$. Such $A$ are called R-recurrent. Let us say that $A$ is
strongly R-positive if the return times of the Markov chain with kernel $P$
have exponential moments of some positive order. Then it is known that strong
R-positivity is equivalent to the property that lowering the value of finitely
many entries of $A$ lowers the spectral radius. This paper gives a short and
largely self-contained proof of this fact.
%
%
\end{abstract}
\vspace{.5cm}

\noi
{\it Acknowledgement.} Work sponsored by grant 16-15238S of the Czech Science
Foundation (GA CR). 



\section{Introduction and main results}

\subsection{R-recurrence}

A nonnegative matrix $A=(A(x,y))_{x,y\in S}$ indexed by a countable set $S$ is
called \emph{irreducible} if for each $x,y\in S$ there exists an $n\geq 1$
such that $A^n(x,y)>0$; it is moreover \emph{aperiodic} if the greatest common
divisor of $\{n\geq 1:A^n(x,x)>0\}$ is one for some, and hence for all $x\in
S$. The classical Perron-Frobenius theorem \cite{Per07,Fro12} says that if $A$
is an irreducible nonnegative matrix indexed by a finite set $S$, then it has
a unique positive eigenfunction. More precisely, there exists a function
$h:S\to(0,\infty)$, which is unique up to scalar multiples, and a unique
constant $c>0$, such that $Ah=ch$. The function $h$ is called the
\emph{Perron-Frobenius eigenfunction} and $c$ the \emph{Perron-Frobenius
  eigenvalue}. We will be interested in generalizations of this theorem to
countably infinite matrices.


Let $A$ be an aperiodic irreducible nonnegative matrix indexed by a countable
set $S$. A simple argument based on superadditivity \cite{Kin63}, shows that
the limit
\be\label{rhodef}
\rho(A):=\lim_{n\to\infty}\big(A^n(x,x)\big)^{1/n}
\ee
exists in $(0,\infty]$ and does not depend on $x\in S$. If $A$ is periodic,
then $\rho(A)$ is defined in the same way except that in (\ref{rhodef}) $n$
ranges only through those integers for which $A^n(x,x)>0$. Because of its
interpretation in the finite case, the quantity $\rho(A)$ is called the
\emph{spectral radius} of $A$. By definition, $A$ is called
\emph{R-recurrent}\footnote{Originally, the letter R was mathematical notation
  for $1/\rho(A)$. For us the `R' in the words R-recurrence, R-positivity
  etc.\ will just be part of the name and not refer to any mathematical
  constant.} if $\rho(A)<\infty$ and
\be\label{Rrec}
\sum_{n=1}^\infty\rho(A)^{-n}A^n(x,x)=\infty
\ee
for some, and hence for all $x\in S$. We observe that a function
$h:S\to(0,\infty)$ is an eigenfunction of $A$ with eigenvalue $c>0$ if and
only if
\be\label{PsimA}
P(x,y):=c^{-1}h(x)^{-1}A(x,y)h(y)\qquad(x,y\in S)
\ee
defines a probability kernel on $S$. As will be shown in Appendix~\ref{A:Rrec}
below, the following theorem follows easily from the work of Vere-Jones
\cite{Ver62,Ver67}.

\bt[R-recurrent matrices]
Let\label{T:Rrec} $A$ be an R-recurrent irreducible nonnegative matrix indexed
by a countable set $S$. Then there exists a function $h:S\to(0,\infty)$, which
is unique up to scalar multiples, and a unique constant $c>0$, such that
(\ref{PsimA}) defines a recurrent probability kernel $P$. Moreover, $c=\rho(A)$.
\et

Since finite matrices are R-recurrent (this is proved in \cite[Sect.~7]{Ver67}
and will also follow from Theorem~\ref{T:strongR} below), the classical
Perron-Frobenius theorem is implied by Theorem~\ref{T:Rrec}. If $S$ is finite,
then $\rho(A)=\sup\{|\la|:\la\in\sig(A)\}$, where $\sig(A)$ denotes the set of
all complex eigenvalues of $A$; by contrast, if $S$ is infinite, then it often
happens that $A$ has positive eigenfunctions with eigenvalues $c>\rho(A)$
\cite{Ver63}. For such eigenfunctions, the probability kernel in (\ref{PsimA})
is transient. Note that in view of this, the term ``spectral radius'' for
$\rho(A)$ is somewhat of a misnomer if $A$ is infinite, but we retain it for
historical reasons. The following theorem, first proved in
\cite[Thm~4.1]{Ver67}, shows that for R-recurrent matrices, there is only one
positive eigenfunction associated with the eigenvalue $\rho(A)$, and all other
positive eigenfunctions (if there are any) have eigenvalues $c>\rho(A)$.

\bt[Positive eigenfunction]
Let\label{T:poseig} $A$ be an R-recurrent irreducible nonnegative matrix indexed
by a countable set $S$. Then there exists a function $h:S\to(0,\infty)$, which
is unique up to scalar multiples, such that such that $Ah=\rho(A)h$. Moreover,
if some function $f:S\to\half$ satisfies $Af\leq\rho(A)f$, then $f=\la h$ for
some $\la\geq 0$.
\et

It should be noted that the approach based on R-recurrence is just one of many
different ways to generalize the Perron-Frobenius theorem to infinite
dimensions. For a more functional analytic approach, see, e.g.,
\cite{KR48,KR50,Kar59,Sch74,Zer87}. The theory of R-recurrence is treated in
the books \cite{Sen81,Woe00} and generalized to uncountable spaces in
\cite{Num84}.

In view of Theorems~\ref{T:Rrec} and \ref{T:poseig}, it is clearly very useful
to know of a given nonnegative matrix that it is R-recurrent. Unfortunately,
it is often not feasible to check this directly from the definition
(\ref{Rrec}), since this requires rather subtle knowledge about the
asymptotics of the powers of $A$, while often it is not even possible to
obtain the spectral radius $\rho(A)$ in closed form. In the next section, we
will show that for a subclass of the R-recurrent matrices, more robust methods
are available.

\subsection{Strong R-positivity}

Let $X=(X_k)_{k\geq 0}$ be a Markov chain with countable state space $S$ and
transition kernel $P$, and let $\sig_x:=\inf\{k>0:X_k=x\}$ denote its first
return time to a point $x\in S$. Let $\P^x$ denote the law of $X$ started in
$X_0=x$ and let $\E^x$ denote expectation with respect to $\P^x$. Recall that
by definition, $x$ is recurrent if $\P^x[\sig_x<\infty]=1$ and $x$ is positive
recurrent if $\E^x[\sig_x]<\infty$. We will say that $x$ is \emph{strongly
  positive recurrent}\footnote{This should be distinguished from the closely
  related, but different concepts of \emph{strong ergodicity} and \emph{strong
    recurrence}, the latter having been defined in \cite{Spi90}.} if
$\E^x[e^{\eps\sig_x}]<\infty$ for some $\eps>0$. It is well-known that
recurrence and positive recurrence are class properties. Kendall \cite{Ken59}
proved that the same is true for strong positive recurrence, i.e., if $P$ is
irreducible, then $\{x\in S:x\mbox{ is strongly positive recurrent}\}$ is
either $S$ or $\emptyset$. For aperiodic chains, he moreover proved
that strong positive recurrence is equivalent to \emph{geometric ergodicity},
in the following sense. In Kendall's formulation, the constant $\eps$ in
point~(iii) was allowed to depend on $x,y$. Vere-Jones \cite{Ver62} showed
that it can be chosen uniformly.

\bp[Geometric ergodicity]
Let\label{P:geom} $P$ be the transition kernel of an irreducible, aperiodic,
positive recurrent Markov chain with countable state space $S$, and let $\pi$
denote its invariant law. Then the following statements are equivalent.
\begin{enumerate}
\item $P$ is strongly positive recurrent.
\item There exist $x\in S$, $\eps>0$, and $M<\infty$ such that
  $|P^n(x,x)-\pi(x)|\leq Me^{-\eps n}$ for all $n\geq 0$.
\item There exist $\eps>0$ and $M_{x,y}<\infty$ such that
  $|P^n(x,y)-\pi(y)|\leq M_{x,y}e^{-\eps n}$ for all $n\geq 0$ and $x,y\in S$.
\end{enumerate}
\ep

An R-recurrent irreducible nonnegative matrix $A$ is called \emph{R-positive}
if the unique recurrent probability kernel $P$ from Theorem~\ref{T:Rrec} is
positive recurrent. We will say that $A$ is \emph{strongly R-positive} if $P$
is strongly positive recurrent. (This is called \emph{geometrically
  R-recurrent} in \cite{Num84}.) An irreducible nonnegative matrix $A$ is
called \emph{R-transient} if it is not R-recurrent and \emph{R-null recurrent}
if it is R-recurrent but not R-positive. We will say that $A$ is \emph{weakly
  R-positive} if it is R-positive but not strongly so. The main aim of the
present paper is to give a short and reasonably self-contained proof of the
following theorem, that gives necessary and sufficient conditions for strong
R-positivity.


\bt[Strong R-positivity]
Let\label{T:strongR} $A$ be an irreducible nonnegative matrix indexed by a
countable set $S$ and assume that $\rho(A)<\infty$. Let $B\leq A$ be another
nonnegative matrix such that $B(x,y)>0$ if and only if $A(x,y)>0$ $(x,y\in
S)$, and $B\neq A$. Then:
\begin{itemize}
\item[{\rm\textbf{(a)}}] If $A$ is strongly R-positive, then $\rho(B)<\rho(A)$.
\item[{\rm\textbf{(b)}}] If $\rho(B)<\rho(A)$ and the set
  $\{(x,y)\in S^2:A(x,y)\neq B(x,y)\}$ is finite, then $A$ is strongly
  R-positive.
\end{itemize}
\et

Theorem~\ref{T:strongR} says that a nonnegative matrix is strongly R-positive
if and only if lowering the value of finitely many entries lowers the spectral
radius. In view of this, to prove strong R-positivity, it suffices to prove
sufficiently sharp upper and lower bounds on the spectral radii of two
nonnegative matrices, which is in general much easier than determining the
exact asymptotics as in (\ref{Rrec}).

For R-transience, a complementary statement holds. The following theorem says
that a nonnegative matrix is R-transitive if and only if it is possible to
increase the value of finitely many entries without increasing the spectral
radius.

\bt[R-transience]
Let\label{T:Rtrans} $A\leq B$ be irreducible nonnegative matrices indexed by a
countable set $S$ and assume that $B\neq A$ and $\rho(B)<\infty$. Then:
\begin{itemize}
\item[{\rm\textbf{(a)}}] If $A$ is R-transient and $\{(x,y)\in S^2:A(x,y)\neq
  B(x,y)\}$ is finite, then $\rho\big(A+\eps(B-A)\big)=\rho(A)$ for some
  $\eps>0$.
\item[{\rm\textbf{(b)}}] If $\rho(B)=\rho(A)$, then $A$ is R-transient.
\end{itemize}
\et

\subsection{Discussion}\label{S:discus}

For nonnegative matrices $A,B$, let us write $A\asymp B$ if there exists a
constant $C\in(0,\infty)$ such that $C^{-1}A\leq B\leq CA$, and let us say
that $A,B$ are \emph{finite modifications} of each other if $A\asymp B$ and
$\{(x,y)\in S^2:A(x,y)\neq B(x,y)\}$ is finite. Theorem~\ref{T:Rtrans}
says that an irreducible nonnegative matrix $A$ is strongly R-positive if and
only if $\rho(B)<\rho(A)$ for some, and hence for all finite modifications $B$
of $A$ such that $B\neq A$. In practice, when the aim is to prove strong
R-positivity for $A$, it is best to choose $B$ as small as possible. One can
go one step further and define
\be
\rho_\infty(A):=\inf\big\{\rho(B):B\leq A,\ B\mbox{ is a finite modification
  of }A\big\}.
\ee
Then Theorem~\ref{T:strongR} implies that $A$ is stongly R-positive if and
only if $\rho_\infty(A)<\rho(A)$. We can describe this condition in words by
saying that under a Gibbs measure with transfer matrix $A$, paths far from the
origin carry less mass, on an exponential scale, than paths near the
origin. The quantity $\rho_\infty(A)$ has been studied in
\cite{MS95,Ign06}. In the latter paper, it is called the \emph{essential
  spectral radius}.

Theorem~\ref{T:strongR} has an interesting history. If $A,B$ are nonnegative
matrices indexed by countable sets $S,T$, respectively, then let us say that
$B$ is a \emph{submatrix} of $A$ if $T\sub S$ and $B_{ij}\in\{0,A_{ij}\}$ for
all $i,j\in T$, i.e., $B$ is obtained by removing some rows and corresponding
columns from $A$ and by replacing some entries of $A$ by zeros. If moreover
$B\neq A$, then we call $B$ a \emph{proper submatrix}.
The following is \cite[Thm~3.15]{GS98}.

\bt[Characterization in terms of submatrices]
An\label{T:submat} irreducible nonnegative matrix $A$ indexed by a countable set
is strongly R-positive if and only if $\rho(B)<\rho(A)$ for all irreducible
proper submatrices $B$ of $A$.
\et

For matrices with values in $\{0,1\}$, this was proved by Salama
\cite{Sal88}. Ruette \cite[Remark~2.1]{Rue03} pointed out mistakes in
Salama's proofs, fixed them, and extended his results. The proof given in
\cite{GS98} differs from the previous proofs, but it is claimed the original
proof can be modified to obtain the general result. If the aim is to prove
strong R-positivity for $A$, then a drawback of Theorem~\ref{T:submat} is that
one has to check that $\rho(B)<\rho(A)$ for \emph{all} irreducible proper
submatrices $B$ of $A$. Closer in spirit to Theorem~\ref{T:strongR} is
\cite[Remark~3.16]{GS98}, which states that an irreducible nonnegative matrix
$A$ is strongly R-positive if and only if lowering the value of a single entry
$A_{ij}$ lowers the spectral radius. Even though this does not immediately
imply the result for finite modifications, it shows that results in the spirit
of Theorem~\ref{T:strongR}, and the methods needed to prove it, are known.

It is interesting that while the theory of R-recurrence originally rose from
the field of probability theory, most of the recent developments come from
ergodic theory, and more specifically from the theory of countable Markov
shifts. A good modern introduction to this field can be found in \cite{Sar15}.
Much of the theory of R-recurrence can be generalized from nonnegative
matrices to the more general Ruelle operators. In this context, strong
R-positivity corresponds to positive recurrence with the spectral gap
property, for which the Discriminant Theorem \cite[Thm~6.7]{Sar15} gives
necessary and sufficient conditions. Theorems~\ref{T:strongR} and
\ref{T:Rtrans} are also close in spirit to \cite[Thms 2.2 and 2.3]{CS09}.



One possible application of Theorem~\ref{T:strongR} in the study of
quasi-stationary laws. The usefulness of R-positivity in the study of
quasi-stationary laws has been noticed long ago \cite[Thm~3.2]{SV66}; see also
\cite[Prop~5.2.10 en 5.2.11]{And91} for the continuous-time case.

Another application of Theorems~\ref{T:strongR} and \ref{T:Rtrans} is in the
study of pinning models. In fact, using these theorems, it is easy to
prove that for pinning models in the localized regime, return times have
exponential moments of some positive order. Moreover, at the critical point
separating the localized and delocalized regimes, the model is either null
recurrent or weakly positive recurrent. These facts have been noticed before,
see \cite[Thm~4.1 and Prop.~4.2]{CGZ06} and \cite[Thm~2.3]{Gia07}.
Note that these references do not explicitly discuss exponential moments of
return times; nevertheless, the claims follow from their formulas.

As demonstrated by the following corollaries, Theorem~\ref{T:strongR} can also
be used to prove strong positive recurrence (which, by
Proposition~\ref{P:geom}, in the aperiodic case is equivalent to geometric
ergodicity).

\bcor[Conditions for strong positive recurrence]
Let\label{C:strong} $P$ be an irreducible probability kernel $P$ on a countable
set $S$ and let $Q\leq P$ be a subprobability kernel such that $Q(x,y)>0$ if
and only if $P(x,y)>0$, and $\{(x,y)\in S^2:Q(x,y)<P(x,y)\}$ is finite.
Then $P$ is strongly positive recurrent if and only if $\rho(Q)<\rho(P)=1$.
\ecor
\bpro
It is easy to see that a recurrent probability kernel $P$ satisfies
$\rho(P)=1$. Moreover, a probability kernel that is strongly positive
recurrent is clearly strongly R-positive. In view of this, the necessity of
the conditions $\rho(Q)<\rho(P)=1$ follows from Theorem\ref{T:strongR}~(a).

Conversely, if $\rho(Q)<\rho(P)$, then Theorem\ref{T:strongR}~(b) shows that
$P$ is strongly R-positive, i.e., there exists a strongly positive recurrent
probability kernel $P'$ and a function $h:S\to(0,\infty)$ such that
$P'(x,y)=\rho(P)^{-1}h(x)^{-1}P(x,y)h(y)$. It follows that $Ph=\rho(P)h$.
Since $\rho(P)=1$, the constant function $h\equiv 1$ also solves this
equation. Since by Theorem~\ref{T:poseig}, solutions of $Ph=\rho(P)h$ are up
to a multiplicative constant unique, we conclude that $P=P'$.
\epro

\bcor[Conditions for strong positive recurrence]
Let\label{C:strong2} $P$ be an irreducible probability kernel $P$ on a
countable set $S$ and let $S'\sub S$ be finite. Then $P$ is strongly positive
recurrent if and only if the following conditions is satisfied:
\begin{enumerate}
\item $\rho(P)=1$
\item There exists a function $f:S\to(0,\infty)$ and $0<\eps<1$ such that
$Pf(x)<\infty$ for all $x\in S'$ and $Pf(x)\leq(1-\eps)f(x)$ for all $x\in
S\beh S'$.
\end{enumerate}
\ecor
\bpro
If (ii) holds, then we can construct a subprobability kernel $Q$ with 
$Q(x,y)>0$ if and only if $P(x,y)>0$ and $Q(x,y)=P(x,y)$ for all $x\in S\beh
S'$, such that $\{(x,y)\in S^2:Q(x,y)<P(x,y)\}$ is finite and
\be
\sum_yQ(x,y)f(y)\leq(1-\eps)f(x)\qquad(x\in S').
\ee
Then $Qf\leq(1-\eps)f$, which is easily seen to imply
$\rho(Q)\leq1-\eps$. Together with condition~(i), by Corollary~\ref{C:strong},
this implies the strong positive recurrence of $P$.

Conversely, if $P$ is strongly positive recurrent, pick some $x_0\in S'$ and
$y_0\in S$ with $P(x_0,y_0)>0$ and define $Q(x_0,y_0):=\ha P(x_0,y_0)$ and
$Q(x,y):=P(x,y)$ for all $(x,y)\neq(x_0,y_0)$. By Corollary~\ref{C:strong},
$\rho(Q)<1$, so by Lemma~\ref{L:exces} below there exists a function
$f:S\to(0,\infty)$ such that $Qf\leq\rho(Q)f$ and hence $Pf(x)\leq\rho(Q)f(x)$
for all $x\neq x_0$, while $Pf(x_0)\leq 2Qf(x_0)<\infty$.
\epro

Corollary~\ref{C:strong2} is similar to a result of Popov \cite{Pop77}, who
proved that the function $f$ in condition~(ii) can be chosen such that $f\geq
1$, and with this extra condition on $f$, condition~(i) can be dropped.





The rest of the paper is dedicated to proofs.
%
%
The organization is as follows. Section~\ref{S:notat} contains
some preliminary definitions and lemmas. Section~\ref{S:logmom} gives a
characterization of forms of R-recurrence in terms of a logarithmic moment
generating function.  Using this, Theorems~\ref{T:strongR} and \ref{T:Rtrans}
are then proved in Sections~\ref{S:strocond} and \ref{S:tracond},
respectively. In Appendix~\ref{A:Rrec}, it is explained how
Theorem~\ref{T:Rrec} follows from the work of
Vere-Jones. Appendix~\ref{A:logmom} contains some general facts about
logarithmic moment generating functions.

\section{Proofs}

\subsection{Excursions away from subgraphs}\label{S:notat}

Given a nonnegative matrix $A$ indexed by a countable set $S$, we define a
directed graph $G=(S,E)$ with vertex set $S$ and set of directed edges $E$
given by $E:=\{(x,y)\in S^2:A(x,y)>0\}$. Alternatively, we denote an edge by
$e=(x,y)$ and call $e^-:=x$ and $e^+:=y$ its \emph{starting vertex} and
\emph{endvertex}, respectively. A \emph{walk} in $G$ is a function
$\om:\{0,\ldots,n\}\to S$ with $n\geq 0$ such that
\be
\vec\om_k:=(\om_{k-1},\om_k)\in E\qquad(1\leq k\leq n).
\ee
We call $\ell_\om:=n\geq 0$ the \emph{length} of $\om$ and we call
$\om^-:=\om_0$ and $\om^+:=\om_n$ its \emph{starting vertex} and
\emph{endvertex}. We can, and sometimes will, naturally identify walks of
length zero and one with vertices and edges, respectively. We let $\Om=\Om(G)$
denote the space of all walks in $G$ and write
\be
\Om^n:=\{\om\in\Om:\ell_\om=n\}
\quand
\Om_{x,y}:=\{\om\in\Om:\om^-=x,\ \om^+=y\}
\ee
and $\Om^n_{x,y}:=\Om^n\cap\Om_{x,y}$.
We observe that
\be\label{Aidef}
A^n(x,y)=\sum_{\om\in\Om^n_{x,y}}\Ai(\om)
\quad\mbox{with}\quad
\Ai(\om):=\prod_{k=1}^{\ell_\om}A(\om_{k-1},\om_k).
\ee
This formula also holds for $n=0$ provided we define the empty product as $:=1$.

If $S'\sub S$ is a subset of vertices, then an \emph{excursion away from} $S'$
is a walk $\om\in\Om$ of length $\ell_\om\geq 1$ such that $\om^\pm\in S'$ and
$\om_k\not\in S'$ for all $0<k<\ell_\om$. We denote the set of all excursions
away from $S'$ by $\hatOm(S')$. We sometimes view a graph as the disjoint
union of its vertex and edge sets, $G=S\cup E$. A \emph{subgraph} of $G$ is
then a set $F\sub G$ such that $e^\pm\in S\cap F$ for all $e\in E\cap F$.
Extending our earlier definition, an \emph{excursion away from} $F$ is an
element $\om\in\hatOm(F\cap S)$ such that moreover $\om\not\in F\cap E$, where
we naturally identify edges with walks of length one. We denote the set of all
excursions away from $F$ by $\hatOm(F)$ and write
$\hatOm_{x,y}(F):=\hatOm(F)\cap\Om_{x,y}$, $\hatOm^n(F):=\hatOm(F)\cap\Om^n$,
etc.

For each subgraph $F$ of $G$ and $x,y\in S\cap F$, we define a \emph{moment
  generating function} $\phi^F_{x,y}$ and \emph{logarithmic moment generating
  function} $\psi^F_{x,y}$ by
\be\label{momgenf}
\phi^F_{x,y}(\la):=\sum_{\om\in\hatOm_{x,y}(F)}\ex{\la\ell_\om}\Ai(\om)
\quand
\psi^F_{x,y}(\la):=\log\phi^F_{x,y}(\la)
\qquad(\la\in\R).
\ee
Here $\phi^F_{x,y}$ and $\psi^F_{x,y}$ may be $\infty$ for some values of
$\la$; in addition, $\psi^F_{x,y}(\la):=-\infty$ if $\phi^F_{x,y}(\la)=0$.
The following lemma lists some elementary properties of $\psi^F_{x,y}$.

\bl[Logarithmic moment generating functions]
Assume\label{L:psiprop} that $A$ is irreducible and $\rho(A)<\infty$.
Let $F$ be a subgraph of $G$, let $x,y\in S\cap F$, and set
\bc\label{lach}
\dis\la_+=\la^F_{x,y,+}&:=&\dis\sup\{\la\in\R:\psi^F_{x,y}(\la)<\infty\},\\[5pt]
\dis\la_\ast=\la^F_{x,y,\ast}&:=&\dis\sup\{\la\in\R:\psi^F_{x,y}(\la)<0\}.
\ec
Then either $\psi^F_{x,y}\equiv-\infty$ or:
\begin{enumerate}
\item $\psi^F_{x,y}$ is convex.
\item $\psi^F_{x,y}$ is lower semi-continuous.
\item $-\infty<\la_\ast<\infty$ and $\la_\ast\leq\la_+\leq\infty$. 
\item $\psi^F_{x,y}$ is infinitely differentiable on $(-\infty,\la_+)$.
\item $\psi^F_{x,y}$ is strictly increasing on $(-\infty,\la_+)$.
\item $\lim_{\la\to\pm\infty}\psi^F_{x,y}(\la)=\pm\infty$.
\end{enumerate}
\el

\bpro
If $\Ai(\om)=0$ for all $\om\in\hatOm_{x,y}(F)$, then
$\psi^F_{x,y}\equiv-\infty$, while otherwise $\psi^F_{x,y}(\la)>-\infty$
for all $\la\in\R$. Clearly, $\psi^F_{x,y}(\la)$ is nondecreasing as
a function of $\la$. Since
\be
\phi^F_{x,y}(\la)
=\sum_{k=1}^\infty\sum_{\om\in\hatOm^k_{x,y}(F)}\ex{\la k}\Ai(\om)
\leq\sum_{k=0}^\infty\sum_{\om\in\Om^k_{x,y}}\ex{\la k}\Ai(\om)
=\sum_{k=0}^\infty\ex{\la k}A^k(x,y),
\ee
which by (\ref{rhodef}) is finite for $\la<-\log\rho(A)$, we see that
$-\infty<-\log\rho(A)\leq\la_+$. Properties (i)--(iv), except for the fact
that $-\infty<\la_\ast<\infty$, now follow from general properties of
logarithmic moment generating functions, see Lemma~\ref{L:logmom} in the
appendix. Property~(vi) follows by monotone convergence and this implies
$-\infty<\la_\ast<\infty$. Since excursions have length $\geq 1$, formula
(\ref{meanvar}) from Lemma~\ref{L:logmom} moreover implies property~(v).
\epro

The following two lemmas allow us to prove properties of $\phi^F_{x,y}$ for
finite subgraphs $F$ by induction on the number of vertices and edges.

\bl[Removal of an edge]\label{L:eremove}
Let $A$ be a nonnegative matrix, let $G=(S,E)$ be its associated graph, and
let $F$ be a subgraph of $G$. Let $e\in F\cap E$ and let $F':=F\beh\{e\}$. Then
\be
\phi^{F'}_{x,y}(\la)
=\left\{\ba{ll}
\phi^F_{x,y}(\la)+e^\la A(x,y)\quad&\mbox{if }e=(x,y),\\[5pt]
\phi^F_{x,y}(\la)\quad&\mbox{otherwise}
\ea\right.\qquad(\la\in\R).
\ee
\el
\bpro
This is immediate from the definition of the moment generating function in
(\ref{momgenf}) and the fact that
\be
\hatOm_{x,y}(F')=\left\{\ba{ll}
\hatOm_{x,y}(F)\cup\{e\}\quad&\mbox{if }e=(x,y),\\[5pt]
\hatOm_{x,y}(F)\quad&\mbox{otherwise,}
\ea\right.
\ee
where we identify $e$ with the walk of length 1 that jumps through $e$.
\epro

\bl[Removal of an isolated vertex]\label{L:xremove}
Let $A$ be a nonnegative matrix, let $G=(S,E)$ be its associated graph, and
let $F$ be a subgraph of $G$. Let $z\in F\cap S$ be a vertex of
$F$. Assume that no edges in $F\cap E$ start or end at $z$ and hence
$F':=F\beh\{z\}$ is a subgraph of $G$. Then
\be\label{xremove}
\phi^{F'}_{x,y}(\la)
=\phi^F_{x,y}(\la)+\sum_{k=0}^\infty
\phi^F_{x,z}(\la)\phi^F_{z,z}(\la)^k\phi^F_{z,y}(\la)
\qquad(x,y\in F'\cap S,\ \la\in\R).
\ee
\el
\bpro
Distinguishing excursions away from $F'$ according to how often they visit the
vertex $z$, we have
\be\ba{r@{\,}l}
\dis\phi^{F'}_{x,y}(\la)
=\sum_{\om_{x,y}}\ex{\la\ell_{\om_{x,y}}}\Ai(\om_{x,y})\\[5pt]
\dis+\sum_{k=0}^\infty
\sum_{\om_{x,z}}\sum_{\om_{z,y}}\sum_{\om^1_{z,z}}\cdots\sum_{\om^k_{z,z}}
&\dis\ex{\la(\ell_{\om_{x,z}}+\ell_{\om_{z,y}}+\ell_{\om^1_{z,z}}
+\cdots+\ell_{\om^k_{z,z}})}\\[-5pt]
&\dis\hspace{2cm}\times
\Ai(\om_{x,z})\Ai(\om_{z,y})\Ai(\om^1_{z,z})\cdots\Ai(\om^k_{z,z}),
\ec
where we sum over $\om_{x,y}\in\hatOm_{x,y}(F)$ etc. Rewriting gives
\be\ba{l}
\dis\phi^{F'}_{x,y}(\la)=\sum_{\om_{x,y}}\ex{\la\ell_{\om_{x,y}}}\Ai(\om_{x,y})\\[5pt]
\dis+\big(\sum_{\om_{x,z}}\ex{\la\ell_{\om_{x,z}}}\Ai(\om_{x,z})\big)
\big(\sum_{\om_{z,y}}\ex{\la\ell_{\om_{z,y}}}\Ai(\om_{z,y})\big)
\sum_{k=0}^\infty\big(\sum_{\om_{z,z}}\ex{\la\ell_{\om_{z,z}}}\Ai(\om_{z,z})\big)^k,
\ec
which is the formula in the lemma.
\epro

\subsection{Excursions away from a single point}\label{S:logmom}

Let $A$ be a nonnegative matrix with index set $S$ and let $G=(S,E)$ be its
associated directed graph. For $z\in S$, we let
\be\label{psidef}
\psi_z:=\psi^{\{z\}}_{z,z},\quad\la_{z,+}:=\la^{\{z\}}_{z,z,+},
\quand\la_{z,\ast}:=\la^{\{z\}}_{z,z,\ast},
\ee
denote the logarithmic moment generating function defined in (\ref{momgenf})
and the constants from (\ref{lach}) for the subgraph $F=\{z\}$ which consists
of the vertex $z$ and no edges. We also write $\hatOm_z:=\hatOm_{z,z}(\{z\})$
for the space of all excursions away from $z$. The following proposition links
forms of R-recurrence to the shape of $\psi_z$.

\bp[Forms of R-recurrence]
Assume\label{P:psichar} that $A$ is irreducible with $\rho(A)<\infty$,
and let $z\in S$ be any reference point. Then
\begin{itemize}
\item[{\rm(a)}] $\la_{z,\ast}=-\log\rho(A)=:\la_\ast$.
\item[{\rm(b)}] One has $\psi_z(\la_\ast)<0$ if $A$ is R-transient and
  $\psi_z(\la_\ast)=0$ if $A$ is R-recurrent.
\item[{\rm(c)}] $A$ is R-positive if and only if the left derivative of $\psi_z$
  at $\la_\ast$ is finite.
\item[{\rm(d)}] $A$ is strongly R-positive if and only if $\la_\ast<\la_{z,+}$.
\end{itemize}
\ep

To prove Proposition~\ref{P:psichar}, we need one preparatory definition and
lemma. Given a nonnegative matrix $A$, for each $\la\in\R$, we define a
\emph{Green's function} $G_\la(x,y)$ by
\be
G_\la(x,y):=\sum_{k=0}^\infty\ex{\la k}A^k(x,y)\qquad(x,y\in S),
\ee
which may be infinite for some values of $\la$. If $A$ is irreducible, then it
is known that \cite[Thm~A]{Ver67}
\be\label{radconv}
G_\la(x,y)<\infty\mbox{ for }\la<-\log\rho(A)
\quand
G_\la(x,y)=\infty\mbox{ for }\la>-\log\rho(A)
\ee
for all $x,y\in S$. The following lemma makes a link between the Green's
function and~$\psi_z$.

\bl[Value on the diagonal]
One\label{L:GPrz} has
\be\label{Gzz}
G_\la(z,z)=\left\{\ba{ll}
(1-e^{\psi_z(\la)})^{-1}\quad&\mbox{if }\psi_z(\la)<0,\\[5pt]
\infty\quad&\mbox{otherwise.}
\ea\right.
\ee
\el
\bpro
Since each $\om\in\Om_{z,z}$ can be written as the
concatenation of $m\geq 0$ excursions $\om^{(i)}\in\hatOm_z$, using the
convention that a product of $m=0$ factors is 1, we see that
\[\ba{l}
\dis\sum_{n=0}^\infty\ex{\la n}A^n(z,z)
=\sum_{n=0}^\infty\sum_{\om\in\Om^n_{z,z}}\ex{\la\ell_\om}\Ai(\om)
=\sum_{\om\in\Om_{z,z}}\ex{\la\ell_\om}\Ai(\om)\\[5pt]
\dis\quad=\sum_{m=0}^\infty\prod_{i=1}^m
\Big(\sum_{\om^{(i)}\in\hatOm_z}\ex{\la\ell_{\om^{(i)}}}\Ai(\om^{(i)})\Big)
=\sum_{m=0}^\infty\ex{m\psi_z(\la)},
\ea\]
which yields (\ref{Gzz}).
\epro

\bpro[of Proposition~\ref{P:psichar}]
Part~(a) follows from formula (\ref{radconv}) and Lemma~\ref{L:GPrz}.
If $A$ is R-transient, then it is immediate from the definition of
R-transience (\ref{Rrec}) that $G_{\la_\ast}(z,z)<\infty$ and hence by
Lemma~\ref{L:GPrz} $\psi_z(\la_\ast)<0$.

On the other hand, if $A$ is R-recurrent, then by Theorem~\ref{T:poseig} there
exists a function $h:S\to(0,\infty)$, which is unique up to scalar multiples,
such that such that $Ah=\rho(A)h$. Setting
\be
P(x,y):=\rho(A)^{-1}h(x)^{-1}A(x,y)h(y)\qquad(x,y\in S)
\ee
now defines a probability kernel. Since $P^n(x,x)=\rho(A)^{-n}A^n(x,x)$,
we see from (\ref{Rrec}) that $\sum_nP^n(x,x)=\infty$ $(x\in S)$, which proves
that $P$ is recurrent. The Markov chain with transition kernel $P$ makes
i.i.d.\ excursions away from $z$ with common law
\be\label{reclaw}
\Pc(\om)=\prod_{k=1}^{\ell_\om}P(\om_{k-1},\om_k)=\rho(A)^{-\ell_\om}\Ai(\om)
=\ex{\la_\ast\ell_\om}\Ai(\om)\qquad(\om\in\hatOm_z).
\ee
In particular, since $P$ is recurrent,
\be
1=\sum_{\om\in\hatOm_z}\Pc(\om)=\sum_{\om\in\hatOm_z}\ex{\la_\ast\ell_\om}\Ai(\om)
=\ex{\psi_z(\la_\ast)}.
\ee
This shows that $\psi_z(\la_\ast)=0$, completing the proof of part~(b).

It follows from Lemmas~\ref{L:logmom} and \ref{L:psilap} in the appendix that
the left derivative of $\psi_z$ at $\la_\ast$ is the mean length of excursions
away from $z$ under the law in (\ref{reclaw}), proving part~(c). Moreover,
\be
\sum_{\om\in\hatOm_z}\Pc(\om)\ex{\eps\ell_\om}=
\sum_{\om\in\hatOm_z}\ex{(\la_\ast+\eps)\ell_\om}\Ai(\om)
=\psi_z(\la_\ast+\eps)
\ee
is finite for $\eps>0$ suffiently small if and only if $\la_\ast<\la_{z,+}$,
proving part~(d).
\epro

\subsection{Characterization of strong R-recurrence}\label{S:strocond}

In this section, we prove Theorem~\ref{T:strongR}. We start with two
preparatory results.

\bl[Excessive functions]
Let\label{L:exces} $A$ be an irreducible nonnegative matrix indexed by a
countable set $S$. Then there exists a function $h:S\to(0,\infty)$ such that
$Ah\leq\rho(A)h$.
\el
\bpro
We can without loss of generality assume that $\rho(A)<\infty$.
If $A$ is R-recurrent, then the statement follows from Theorem~\ref{T:poseig}.
If $A$ is R-transient, then $G_{\la_\ast}(x,x)<\infty$ for each $x\in S$ by
Proposition~\ref{P:psichar}~(b) and Lemma~\ref{L:GPrz}. A simple argument based
on irreducibility shows that $G_{\la_\ast}(x,y)<\infty$ for each $x,y\in
S$. Since
\be
AG_{\la_\ast}(x,z)=\sum_{k=0}^\infty\ex{\la_\ast k}A^{k+1}(x,z)
=\ex{-\la_\ast}G_{\la_\ast}(x,z)-1_{\{x=z\}}=\rho(A)G_{\la_\ast}(x,z)-1_{\{x=z\}},
\ee
setting $h(x):=G_{\la_\ast}(x,z)$ $(x\in S)$, where $z\in S$ is any reference
point, now proves the claim.
\epro

\bp[Exponential moments of excursions]\label{P:subMplus}
Let $P$ be an irreducible subprobability kernel on a countable set $S$. Let
$G=(S,E)$ be the graph associated with $P$ and for any subgraph $F\sub G$, let
$\la^F_{x,y,+}$ be defined in terms of $P$ as in (\ref{lach}). Then, if
\be\label{exmco}
\la^F_{x,y,+}>0\mbox{ for all }x,y\in F\cap S
\ee
holds for some finite nonempty subgraph $F$ of $G$, it holds for all such
subgraphs.
\ep
\bpro
We need to show that if $F,F'$ are finite nonempty subgraphs of $G$, then
(\ref{exmco}) holds for $F$ if and only if it holds for $F'$. It suffices to
consider only the following two cases: I.\ $F'=F\beh\{e\}$ where $e\in
F\cap E$ is some edge in $F$, and II.\ $F'=F\beh\{z\}$ where $z\in F\cap S$ is
an isolated vertex in $F$. By Lemma~\ref{L:eremove}, removing an edge does not
change the value of $\la^F_{x,y,+}$ for any $x,y\in F\cap S$, so case~I is
easy.

In case~II, we first prove that if (\ref{exmco}) holds for $F$, then it also
holds for $F'$. We distinguish two subcases: II.a: there exists an
$\om\in\hatOm_{x,y}(F')$ 
that passes through $z$, and II.b: no such $\om$ exists. In case~II.b,
\be\label{laplu2}
\la^{F'}_{x,y,+}=\la^F_{x,y,+},
\ee
so this case is trivial. In case~II.a, Lemma~\ref{L:xremove} tells us that
\be\label{laplu}
\la^{F'}_{x,y,+}=\la^F_{x,y,+}\wedge\la^F_{x,z,+}\wedge\la^F_{z,y,+}
\wedge\la^F_{z,z,\ast}.
\ee
We claim that
\be\label{astplus}
\la^F_{z,z,\ast}>0\quad\mbox{if and only if}\quad\la^F_{z,z,+}>0.
\ee
To see this, we observe that $\sum_{\om\in\hatOm^F_{z,z}}\Pc(\om)$ is the
probability that the Markov chain with transition kernel $P$ started in
$z$ returns to $z$ before visiting any point of $F'$ or being killed.
Since there exists an $\om\in\hatOm_{x,y}(F')$ 
that passes through $z$, this probability is $<1$ and hence
\[
\psi^F_{z,z}(0)=\log\Big(\sum_{\om\in\hatOm^F_{z,z}}\Pc(\om)\Big)<0.
\]
By Lemma~\ref{L:psiprop}~(i) and (ii), $\psi^F_{z,z}$ is continuous on
$(-\infty,\la^F_{z,z,+}]$, so if $\psi^F_{z,z}(\la)<\infty$ for some $\la>0$
then also $\psi^F_{z,z}(\la)<0$ for some $\la>0$, proving (\ref{astplus}).
Combining (\ref{laplu}) and (\ref{astplus}), we see that if (\ref{exmco})
holds for $F$, then it also holds for $F'$.

We next show that if (\ref{exmco}) does not hold for $F$, then neither does it
for $F'$. We consider four cases: (i) $\la^F_{x,y,+}\leq 0$ for some $x,y\in
F'\cap S$, (ii) $\la^F_{x,z,+}\leq 0$ for some $x\in F'\cap S$, (iii)
$\la^F_{z,y,+}\leq 0$ for some $y\in F'\cap S$, and (iv) $\la^F_{z,z,+}\leq
0$. In case~(i), formulas (\ref{laplu2}) and (\ref{laplu}) immediately show
that $\la^{F'}_{x,y}\leq 0$.  In case~(ii), by irreducibility, we can find
some $y\in F'$ and $\om\in\hatOm_{z,y}(F)$; 
then we are in case~II.a and (\ref{laplu}) implies that $\la^{F'}_{x,y}\leq
0$. Case~(iii) is similar to case~(ii). In case~(iv), finally, by
irreducibility we can find $x,y\in F'$ and an $\om\in\hatOm_{x,y}(F')$
that passes through $z$, so (\ref{laplu}) and (\ref{astplus})
imply that $\la^{F'}_{x,y}\leq 0$.
\epro

\bpro[of Theorem~\ref{T:strongR}]
Pick any reference vertex $z\in S$ and let $\psi_z$ be the logarithmic moment
generating function of $A$, as defined in (\ref{momgenf}) and (\ref{psidef}).
Let $\la_{z,+}$ and $\la_{z,\ast}$ be defined as in (\ref{lach}) and
(\ref{psidef}) and recall from Proposition~\ref{P:psichar}~(a) that
$\la_{z,\ast}=\la_\ast:=-\log\rho(A)$. Define $\psi'_z$,
$\la'_{z,+}$, and $\la'_\ast$ in the same way for~$B$.

It is immediately clear from the definition of $\psi_z$ and irreducibility
that $B\neq A$ implies that $\psi'_z(\la)<\psi_z(\la)$ for all $\la$ such
that $\psi'_z(\la)<\infty$. Since $A$ is strongly R-positive,
Proposition~\ref{P:psichar}~(d) implies that $\la_\ast<\la_{z,+}$. It follows
that $\psi'_z(\la_\ast)<\psi_z(\la_\ast)=0$ while
$\la'_{z,+}\geq\la_{z,+}>\la_\ast$, so by the continuity of $\psi'_z$ on
$(-\infty,\la'_{z,+}]$ we have
\[
-\log\rho(B)=\la'_\ast=\sup\{\la\in\R:\psi'_z(\la)<0\}>\la_\ast
=-\log\rho(A),
\]
which shows that $\rho(B)<\rho(A)$.\med

To prove part~(b), we will show that if $\{(x,y)\in S^2:A(x,y)\neq B(x,y)\}$
is finite and $A$ is not strongly R-positive, then $\rho(B)=\rho(A)$. By
Lemma~\ref{L:exces}, there exists a function $h:S\to(0,\infty)$ such that
$Ah\leq\rho(A)h$. We use this function to define subprobability kernels $P$
and $P'$ by
\be\left.\ba{r@{\,}c@{\,}l}\label{PdefA}
\dis P(x,y)&:=&\dis\rho(A)^{-1}h(x)^{-1}A(x,y)h(y),\\[5pt]
\dis P'(x,y)&:=&\dis\rho(A)^{-1}h(x)^{-1}B(x,y)h(y),
\ea\right\}\quad(x,y\in S).
\ee
Since $P^n(x,x)=\rho(A)^{-n}A^n(x,x)$, we see that $\rho(P)=1$,
and likewise $\rho(P')=\rho(A)^{-1}\rho(B)$. Thus, to prove that
$\rho(B)=\rho(A)$, it suffices to prove that $\rho(P')=1$.

Fix any reference point $z\in S$ and from now on, let $\psi_z$ denote the
logarithmic moment generating function of $P$ (and not of $A$ as before), let
$\la_{z,+}$ and $\la_{z,\ast}$ be as in (\ref{lach}), and let $\psi'_z$,
$\la'_{z,+}$, and $\la'_{z,\ast}$ be the same objects defined for $P'$. By
Proposition~\ref{P:psichar}~(a), $\la_{z,\ast}=\la_\ast:=-\log\rho(P)=0$ and
$\la'_{z,\ast}=\la'_\ast:=-\log\rho(P')$, so we need to show that
$\la'_\ast=0$.

Let us say that two nonnegative matrices are \emph{equivalent} if they are
related as in (\ref{PsimA}) for some function $h:S\to(0,\infty)$ and constant
$c>0$. Then, in the light of Theorem~\ref{T:Rrec}, a nonnegative matrix $A$ is
strongly R-positive if and only if it is equivalent to some (necessarily
unique) strongly positive recurrent probability kernel. Since the 
subprobability kernel $P$ from (\ref{PdefA}) is equivalent to $A$, and 
by assumption, $A$ is not strongly R-positive, it follows that also $P$ is 
not strongly R-positive.
Therefore by Proposition~\ref{P:psichar}~(c), $\la_{z,+}=\la_\ast=0$. Since
$P'$ is a subprobability kernel, we have $\psi'_z(0)\leq 0$. Thus, to prove
that $\la'_\ast=0$, it suffices to show that $\psi'_z(\la)=\infty$ for all
$\la>0$ or equivalently $\la'_{z,+}\leq 0$.

Let $F$ be any finite subgraph of $G$ that contains all edges
$(x,y)$ where $P'(x,y)<P(x,y)$. Let $\phi^F_{x,y}$ and $\psi^F_{x,y}$ be
defined as in (\ref{momgenf}) and $\la^F_{x,y,+}$ as in (\ref{lach}). Since
each excursion away from $F$ has the same weight under $\Pc$ and $\Pc'$, it
does not matter whether we use $P$ or $P'$ to define these
quantities. Applying Proposition~\ref{P:subMplus} to the subgraphs $F$ and
$F':=\{z\}$, we see that
\[
\la_{z,+}\leq 0
\quad\desd\quad
\la^F_{x,y,+}\leq 0\mbox{ for some }x,y\in F\cap S
\quad\desd\quad
\la'_{z,+}\leq 0.
\]
In particular, the fact that $\la_{z,+}=0$ implies $\la'_{z,+}\leq 0$,
as required.
\epro

\subsection{Characterization of R-transience}\label{S:tracond}

In this section we prove Theorem~\ref{T:Rtrans}. We need one preparatory
lemma.

\bl[Strictly excessive functions]
Let\label{L:exces2} $A$ be an R-transient irreducible nonnegative matrix
indexed by a countable set $S$, and let $S'\sub S$ be finite. Then there
exists a function $h:S\to(0,\infty)$ such that $Ah\leq\rho(A)h$ and
$Ah<\rho(A)h$ on $S'$.
\el
\bpro
In the proof of Lemma~\ref{L:exces}, we have seen that setting
$h_z(x):=G_{\la_\ast}(x,z)$ defines a function such that
$Ah_z(x)=\rho(A)h_z(x)-1_{\{x=z\}}$. In view of this, the function
$h:=\sum_{z\in S'}h_z$ has all the required properties.
\epro

\bpro[of Theorem~\ref{T:Rtrans}]
We start by proving part~(a). Let $E':=\{(x,y)\in S^2:A(x,y)\neq B(x,y)\}$.
We will show that there exists a matrix $C\geq A$ with $C(x,y)>A(x,y)$ for all
$(x,y)\in E'$ and $\rho(C)\leq\rho(A)$. Since $A\leq A+\eps(B-A)\leq C$ for
$\eps>0$ small enough, the claim then follows.

Let $S':=\{x\in S:(x,y)\in E'\mbox{ for some }y\in S\}$. By R-transience and
Lemma~\ref{L:exces2}, there exists a function $h:S\to(0,\infty)$ such that
$Ah\leq\rho(A)h$ and $Ah<\rho(A)h$ on $S'$. It follows that
\be
P(x,y):=\rho(A)^{-1}h(x)^{-1}A(x,y)h(y)\qquad(x,y\in S)
\ee
defines a subprobability kernel such that $\sum_yP(x,y)<1$ for all $x\in S'$.
Using this, we can construct a probability kernel $Q\geq P$ such that
$Q(x,y)>P(x,y)$ for all $(x,y)\in E'$. Setting
\be
C(x,y):=\rho(A)h(x)Q(x,y)h(y)^{-1}\qquad(x,y\in S)
\ee
then defines a nonnegative matrix $C\geq A$ with $C(x,y)>A(x,y)$ for all
$(x,y)\in E'$. Since $C^n(x,x):=\rho(A)^nQ^n(x,x)$ and $Q$ is a probability
kernel, we see from (\ref{rhodef}) that $\rho(C)\leq\rho(A)$.


To prove also part~(b), fix a reference point $z\in S$ and let
$\psi_z,\la_{z,+}$, and $\la_{z,\ast}$ be defined in terms of $A$ as in 
(\ref{momgenf}), (\ref{lach}), and (\ref{psidef}). Let $\psi'_z,\la'_{z,+}$, and
$\la'_{z,\ast}$ be the same objects defined in terms of $B$.
By Proposition~\ref{P:psichar}
$\la_{z,\ast}=-\log\rho(A)=-\log\rho(B)=\la'_{z,\ast}$.
The definition of $\psi_z$ and irreducibility imply that
$\psi_z(\la)<\psi'_z(\la)$ for all $\la$ such that $\psi'_z(\la)<\infty$, i.e.,
for $\la\leq\la'_{z,+}$. In particular, this
applies at $\la_\ast=\la'_\ast$ so we see that
$\psi_z(\la_\ast)<\psi'_z(\la'_\ast)\leq 0$. By
Proposition~\ref{P:psichar}~(b), it follows that $A$ is R-transient.
\epro

\appendix

\section{Appendix}

\subsection{R-recurrence}\label{A:Rrec}

\bpro[of Theorem~\ref{T:Rrec}]
By Theorem~\ref{T:poseig}, there exists a function $h:S\to(0,\infty)$, which
is unique up to scalar multiples, such that $Ah=\rho(A)h$. Setting
\be
P(x,y):=\rho(A)^{-1}h(x)^{-1}A(x,y)h(y)\qquad(x,y\in S)
\ee
defines a probability kernel on $S$. Since
$P^n(x,x):=\rho(A)^{-1}A^n(x,x)$, we see from (\ref{Rrec}) that
$\sum_{n=1}^\infty P^n(x,x)=\infty$, proving that $P$ is recurrent.

Conversely, assume that for some function $h:S\to(0,\infty)$ and constant
$c>0$, formula (\ref{PsimA}) defines a recurrent probability kernel. Since $P$
is a probability kernel $Ah=ch$. Since $P^n(x,x):=c^{-n}A^n(x,x)$, it follows
that $\rho(P)=c^{-1}\rho(A)$. Since $P$ is a recurrent probability kernel, for
any $x\in S$,
\be
\sum_{k=1}^\infty e^{\la k}P^k(x,x)<\infty
\quad\desd\quad\la<0
\ee
which shows that $\rho(P)=1$ and hence $c=\rho(A)$. Thus $Ah=\rho(A)h$ and
Theorem~\ref{T:poseig} tells us that $h$ is uniquely determined up to scalar
multiples.
\epro

\subsection{Logarithmic moment generating functions}\label{A:logmom}

Let $\mu$ be a nonzero measure on $\R$.
We define the \emph{logarithmic moment generating function} of $\mu$ as
\be
\psi(\la):=\log\int_R\mu(\di x)e^{\la x}\qquad(\la\in\R),
\ee
where $\log\infty:=\infty$. We write
\be
\Di_\psi:=\{\la\in\R:\psi(\la)<\infty\}
\quand
\Ui_\psi:={\rm int}(\Di_\psi).
\ee

\bl[Logarithmic moment generating functions]\label{L:logmom}
Let $\mu$ be a non\-zero measure on $\R$ and let $\psi$ be its logarithmic
moment generating function. Then
\begin{itemize}
\item[{\rm(i)}] $\psi$ is convex.
\item[{\rm(ii)}] $\psi$ is lower semi-continuous.
\item[{\rm(iii)}] $\psi$ is infinitely differentiable on $\Ui_\psi$.
\end{itemize}
Moreover, for each $\la\in\Di_\psi$, setting
\be
\mu_\la(\di x):=\ex{\la x-\psi(\la)}\mu(\di x)
\ee
defines a probability measure on $\R$ with
\be\left.\ba{rr@{\,}c@{\,}l}\label{meanvar}
{\rm(i)}&\dis\dif{\la}\psi(\la)&=&\dis\langle\mu_\la\rangle,\\[5pt]
{\rm(ii)}&\dis\diff{\la}\psi(\la)&=&\dis\var(\mu_\la),
\ea\right\}\quad(\la\in\Ui_\psi),
\ee
where $\langle\mu_\la\rangle$ and $\var(\mu_\la)$ denote the mean and variance
of $\mu_\la$, respectively.
\el
\bpro
Set
\be\label{Phidef}
\Phi(\la):=\int_\R\mu(\di x)e^{\la x}\qquad(\la\in\R),
\ee
so that $\psi(\la)=\log\Phi(\la)$. We claim that $\la\mapsto\Phi(\la)$ is
infinitely differentiable on $\Ui_\psi$ and
\be\label{difPhi}
\big(\dif{\la}\big)^n\Phi(\la)=\int x^ne^{\la x}\mu(\di x)
\qquad(\la\in\Ui_\psi).
\ee
To justify this, we must show that the interchanging of
differentiation and integral is allowed. We observe that
\be\label{difint}
\dif{\la}\int x^ne^{\la x}\mu(\di x)
=\lim_{\eps\to 0}\int x^n\eps^{-1}(e^{(\la+\eps)x}-e^{\la x})\mu(\di x),
\ee
By our assumption that $\la\in\Ui_\psi$, we can choose $\de>0$ such that
\be\label{2de}
\int_\R\mu(\di x)\big[e^{(\la-2\de)x}+e^{(\la+2\de)x}\big]<\infty.
\ee
Since for any $-\de<\eps<\de$ with $\eps\neq 0$,
\bc
\dis|x|^n\eps^{-1}\big|e^{(\la+\eps)x}-e^{\la x}\big|
&=&\dis|x|^n\Big|\eps^{-1}\!\int_\la^{\la+\eps}\!xe^{\kappa x}\,\di\kappa\Big|\\[8pt]
&\leq&\dis|x|^{n+1}\big[e^{(\la-\de)x}+e^{(\la+\de)x}\big]\qquad(x\in\R),
\ec
which is integrable by (\ref{2de}), we may use dominated
convergence in (\ref{difint}) to interchange the limit and integral.

Since
\be
\int_\R\mu_\la(\di x)=\frac{1}{\Phi(\la)}\ex{\la x}\mu(\di x)=1
\qquad(\la\in\Di_\psi),
\ee
we see that $\mu_\la$ is a probability measure for each $\la\in\Di_\psi$.
Formula (\ref{difPhi}) implies that for each $\la\in\Ui_\psi$
\be\ba{rr@{\,}c@{\,}l}\label{diflog}
{\rm(i)}&\dis\dif{\la}\log\Phi(\la)&=&\dis\dif{\la}\log\int e^{\la x}\mu(\di x)
=\frac{\int xe^{\la x}\mu(\di x)}{\int e^{\la x}\mu(\di x)}
=\li \mu_\la\re,\\[15pt]
{\rm(ii)}&\dis\diff{\la}\log\Phi(\la)&=&\dis\frac{\int x^2e^{\la x}\mu(\di x)
-(\Phi(\la)\int xe^{\la x}\mu(\di x))^2}{\Phi(\la)^2}\\[10pt]
&&=&\dis\int x^2\mu_\la(\di x)
-\Big(\int x\mu_\la(\di x)\Big)^2=\var(\mu_\la),
\ea\ee
proving (\ref{meanvar}).

In particular, if $\mu$ is a compactly supported finite measure, then
$\Ui_\psi=\R$ so these formulas prove that $\psi$ is convex and continuous.
For locally finite $\mu$, we may find compactly supported finite $\mu_n$ such
that $\mu_n\up\mu$ and hence the associated logarithmic moment generating
functions satisfy $\psi_n\up\psi$. Since the $\psi_n$ are convex and
continuous, $\psi$ must be convex and l.s.c. If $\mu$ is not locally finite,
then $\psi(\la)=\infty$ for all $\la\in\R$ so there is nothing to prove.
\epro

The next lemma says that formula (\ref{meanvar})~(i) holds more generally for
$\la\in\Di_\psi$, when the derivative is appropriately interpreted as a
one-sided derivative or limit of derivatives for $\la\in\Ui_\psi$.

\bl[One-sided derivative]\label{L:psilap}
Let $\mu$ be a non\-zero measure on $\R$ and let $\psi$ be its logarithmic
moment generating function. Assume that $\Ui_\psi$ is nonempty,
$\la_+:=\sup\Di_\psi<\infty$, and $\psi(\la_+)<\infty$. Then
\be\label{leftder}
\lim_{\la\up\la_+}\dif{\la}\psi(\la)
=\lim_{\eps\down 0}\eps^{-1}\big(\psi(\la_+)-\psi(\la_+-\eps)\big)
=\li\mu_{\la_+}\re.
\ee
\el
\bpro
Since $\ga\mapsto xe^{x\ga}$ is nondecreasing for each $x\geq 0$, we see that
\be
\eps^{-1}(e^{\la x}-e^{(\la-\eps)x})=\int_{\la-\eps}^\la xe^{x\ga}\di\ga
\ \Big\up\ xe^{\la x}\quad\mbox{as }\eps\down 0\quad(x\geq 0),
\ee
so by monotone convergence, using notation as in (\ref{Phidef}), we see that
\be
\eps^{-1}\big(\Phi(\la_+)-\Phi(\la_+-\eps)\big)
=\int\eps^{-1}(e^{\la_+ x}-e^{(\la_+-\eps)x})\mu(\di x)
\;\Big\up\;\int xe^{\la_+ x}\mu(\di x)\quad\mbox{as }\eps\down 0.
\ee
By assumption $\psi(\la_+)<\infty$ so $\Phi(\la_+)<\infty$, which implies that
$\mu_{\la_+}$ is well-defined, and the second equality in (\ref{leftder}) now
follows as in (\ref{diflog}). Since $\psi$ is convex, this also implies the
first equality in (\ref{leftder}).
\epro

\subsection*{Acknowledgement}

I was motivated to prove Theorem~\ref{T:strongR} based on a suggestion of
Rongfeng Sun, who conjectured a version of this result in joint work with
Gregory Maillard and Nicolas P\'etr\'elis. Francesco Caravenna helped me a lot
orienting myself in the literature on pinning models. I thank Jan Seidler for
pointing out the reference \cite{Ign06}. I am indebted to Sergey Savchenko who
pointed out the references \cite{CS09,GS98,Rue03,Sal88}.

\end{document}